\documentstyle[amscd]{amsart}
\numberwithin{equation}{section}
\theoremstyle{plain}
\newtheorem{thm}{Theorem}[section]

\newtheorem{lem}[thm]{Lemma}
\newtheorem{prop}[thm]{Proposition}
\begin{document}
\title{Cuntz-Krieger algebras associated with Hilbert $C^*$-quad modules 
of commuting matrices
}
\author{Kengo Matsumoto}
\address{ 
Department of Mathematics, 
Joetsu University of Education,
Joetsu 943-8512, Japan}
\email{kengo@@juen.ac.jp}
\maketitle
\begin{abstract}
Let
${\cal O}_{{\cal H}^{A,B}_\kappa}$
be the $C^*$-algebra
associated with
the Hilbert $C^*$-quad module arising from 
 commuting matrices $A,B$ with entries in $\{0,1\}$.
We will show that if the associated tiling space $X_{A,B}^\kappa$
is transitive, the $C^*$-algebra ${\cal O}_{{\cal H}^{A,B}_\kappa}$
is simple and purely infinite.
In particulr, for two positive integers $N,M$, 
the $K$-groups of the simple purely infinite 
$C^*$-algebra ${\cal O}_{{\cal H}^{[N],[M]}_\kappa}$ 
are computed by using the Euclidean algorithm. 
\end{abstract}


\def\Zp{{ {\Bbb Z}_+ }}
\def\CTDS{{ ({\cal A}, \rho, \eta, \Sigma^\rho, \Sigma^\eta, \kappa)}}
\def\OL{{{\cal O}_{{\frak L}}}}
\def\FL{{{\cal F}_{{\frak L}}}}
\def\FKL{{{\cal F}_k^{l}}}
\def\FHK{{{\cal F}_{{\cal H}_\kappa}}}
\def\FU{{{\cal F}_{{\cal H}_\kappa}^{uni}}}
\def\FHLI{{{\cal F}_{[{\frak L}]} }}
\def\FRHO{{ {\cal F}_\rho}}
\def\DRHO{{{{\cal D}_{\rho}}}}
\def\ORHO{{{\cal O}_\rho}}
\def\OETA{{{\cal O}_\eta}}
\def\ORE{{{\cal O}^{\kappa}_{\rho,\eta}}}
\def\OREK{{{\cal O}_{(\rho,\eta;\kappa)}}}
\def\OHK{{{\cal O}_{{\cal H}_\kappa^{\rho,\eta}}}}
\def\OU{{{\cal O}_{{\cal H}_\kappa}^{uni}}}
\def\OAK{{{\cal O}_{\A,\kappa}}}
\def\OHRE{{{\widehat{\cal O}}^{\kappa}_{\rho,\eta}}}
\def\FRE{{{\cal F}_{\rho,\eta}}}
\def\DRE{{{\cal D}_{\rho,\eta}}}
\def\M{{{\cal M}}}
\def\N{{{\cal N}}}
\def\HK{{{\cal H}_\kappa^{\rho,\eta}}}
\def\EK{{E_\kappa}}
\def\SK{{\Sigma_\kappa}}
\def\K{{{\cal K}}}
\def\A{{{\cal A}}}
\def\B{{{\cal B}}}
\def\BE{{{\cal B}_\eta}}
\def\BR{{{\cal B}_\rho}}
\def\BK{{{\cal B}_\kappa}}
\def\BU{{{\cal B}_\kappa^{uni}}}
\def\Aut{{{\operatorname{Aut}}}}
\def\End{{{\operatorname{End}}}}
\def\Hom{{{\operatorname{Hom}}}}
\def\Ker{{{\operatorname{Ker}}}}
\def\ker{{{\operatorname{ker}}}}
\def\Coker{{{\operatorname{Coker}}}}
\def\id{{{\operatorname{id}}}}
\def\dim{{{\operatorname{dim}}}}
\def\min{{{\operatorname{min}}}}
\def\exp{{{\operatorname{exp}}}}
\def\supp{{{\operatorname{supp}}}}
\def\Proj{{{\operatorname{Proj}}}}
\def\Im{{{\operatorname{Im}}}}
\def\span{{{span}}}
\def\OA{{{\cal O}_A}}
\def\OB{{{\cal O}_B}}
\def\T{{  {\cal T}_{{\K}^{\M}_{\N}} }}
\def\TK{{ {\cal T}_{{\cal H}_\kappa} }}
\section{Introduction}
 
 In \cite{MaCrelle}, the author has introduced a notion of 
 $C^*$-symbolic dynamical system, 
 which is a generalization of a finite labeled graph,
  a $\lambda$-graph system and an automorphism of a unital $C^*$-algebra.
It is denoted by $(\A, \rho,\Sigma)$
and consists of a finite family $\{ \rho_{\alpha} \}_{\alpha \in \Sigma}$ 
of endomorphisms of a unital $C^*$-algebra 
${\cal A}$ 
such that 
$\rho_\alpha(Z_\A) \subset Z_\A, \alpha \in \Sigma$
and
$\sum_{\alpha\in \Sigma}\rho_\alpha(1) \ge 1$
where
$Z_\A$ denotes the center of $\A$.
It provides a subshift 
$\Lambda_\rho$  
over $\Sigma$
and a Hilbert $C^*$-bimodule 
${\cal H}_{\cal A}^{\rho}$
over $\A$
which gives rise to a  $C^*$-algebra
${\cal O}_\rho$ as a Cuntz-Pimsner algebra 
(\cite{MaCrelle}, cf. \cite{KPW}, \cite{Pim}).
In \cite{MaPre2011} and \cite{MaPre2011No2}, 
the author has extended the notion of
$C^*$-symbolic dynamical system to 
$C^*$-textile  dynamical system
which is a higher dimensional analogue of $C^*$-symbolic dynamical system.
The $C^*$-textile  dynamical system
$({\cal A}, \rho, \eta, \Sigma^\rho, \Sigma^\eta, \kappa)$
 consists of two 
$C^*$-symbolic dynamical systems 
$({\cal A}, \rho,  {\Sigma^\rho})$
and 
$({\cal A}, \eta, {\Sigma^\eta})$
with a common unital $C^*$-algebra $\A$ and a  commutation relation
between the endomorphisms $\rho$ and $\eta$ through a map $\kappa$ 
stated below.
Set 
\begin{align*}
\Sigma^{\rho \eta} & = \{ (\alpha, b ) \in \Sigma^\rho \times \Sigma^\eta \mid 
\eta_b \circ \rho_\alpha  \ne 0 \},\\
\Sigma^{\eta \rho} & = \{ (a, \beta) \in \Sigma^\eta \times \Sigma^\rho \mid 
\rho_\beta \circ \eta_a  \ne 0 \}.
\end{align*}
We assume that there exists a bijection
$
\kappa : \Sigma^{\rho \eta} \longrightarrow \Sigma^{\eta\rho},
$
which we fix and call a specification.
Then the required commutation relations are 
\begin{equation}
\eta_b \circ \rho_\alpha = \rho_\beta \circ \eta_a
\qquad
\text{ if } 
\quad 
\kappa(\alpha, b) = (a,\beta). \label{eqn:kappa}
\end{equation}
A $C^*$-textile dynamical system
 provides  a two-dimensional subshift and 
a multi structure of Hilbert $C^*$-bimodule that has  multi right actions and 
multi left actions and multi inner products.
Such  a multi structure of Hilbert $C^*$-bimodule is called 
a Hilbert $C^*$-quad module,
 denoted  by
 $\HK$.
In \cite{MaPre2011No2},
the author has introduced 
a $C^*$-algebra associated with the Hilbert $C^*$-quad module.
The  $C^*$-algebra $\OHK$
has been constructed in a concrete way from 
the  structure of the Hilbert $C^*$-quad module
$\HK$ by a two-dimensional analogue of 
Pimsner's construction of $C^*$-algebras from Hilbert $C^*$-bimodules.
It is generated by the quotient images 
of creation operators on two-dimensional analogue of Fock Hilbert 
module by module maps of compact operators. 
As a result,
the $C^*$-algebra has been proved to have a universal property 
subject to certain operator relations 
of generators encoded by structure of the Hilbert $C^*$-quad module (\cite{MaPre2011No2}).

Let
$A,B$ be two $N \times N$ matrices with entries in 
nonnegative integers.
We assume that both $A$ and $B$ are essential,
which means that they have no rows or columns identically to zero vector.
They yield directed graphs
$G_A = (V, E_A)$
and 
$G_B = (V, E_B)$
with a common vertex set
$V = \{ v_1,\dots, v_N \}$ and 
edge sets $E_A$ and 
$E_B$ respectively,
where
the edge set 
$E_A$ consist of $A(i,j)$-edges from the vertex 
$v_i$ to the vertex $v_j$
and 
$E_B$ consist of $B(i,j)$-edges from the vertex 
$v_i$ to the vertex $v_j$.
We then have two
$C^*$-symbolic dynamical systems
$(\A_N, \rho^A, E_A)$
and
$(\A_N, \rho^B, E_B)$
with
$\A_N = {\Bbb C}^N$.
Denote by $s(e), r(e)$
the source vertex and the range vertex of an edge $e$.
Put
\begin{align*}
\Sigma^{AB}
& =\{(\alpha,b)\in E_A \times E_B 
\mid  r(\alpha) = s(b)  \}, \\ 
\Sigma^{BA}
& =\{(a,\beta)\in E_B \times E_A 
\mid  r(a) = s(\beta) \}.
\end{align*}
Assume that the commutation relation
\begin{equation}
AB = BA  \label{eqn:abba}
\end{equation} 
holds.
We may take a bijection 
$\kappa: \Sigma^{AB} \longrightarrow \Sigma^{BA}$
such that $s(\alpha) = s(a), r(b) = r(\beta)$
for $\kappa(\alpha,b) = (a,\beta)$,
which we fix.
This situation is called an LR-textile system introduced by Nasu
(\cite{NaMemoir}).
We then have a $C^*$-textile dynamical system
(see \cite{MaPre2011No2})
\begin{equation*}
(\A_N, \rho^A,\rho^B, E_A, E_B, \kappa).
\end{equation*}
Let us denote  by
${\cal H}_\kappa^{A,B}$ the associated Hilbert $C^*$-quad module
defined in \cite{MaPre2011No2}.
We set
\begin{equation}
E_\kappa
=\{ (\alpha, b, a,\beta) \in E_A \times E_B \times E_B \times E_A | \kappa(\alpha,b) = (a, \beta) \}. \label{eqn:EK}
\end{equation}
Each element of $E_\kappa$ is called a tile.
Let
$X_{A,B}^\kappa \subset (E_\kappa)^{{\Bbb Z}^2}$ 
be the two-dimensional subshift of
the Wang tilings of $E_\kappa$ (cf. \cite{Wang}). 
It consists of the two-dimensional configurations
$x : {\Bbb Z}^2 \longrightarrow E_\kappa$
compatible to their boundary edges on each tile,
and is called the subshift of the tiling space 
for the specification
$\kappa: \Sigma^{AB} \longrightarrow \Sigma^{BA}$.
We say that 
$X_{A,B}^\kappa$ is transitive
if for two tiles
$\omega, \omega' \in E_\kappa$,
there exists 
$(\omega_{i,j})_{(i,j) \in {\Bbb Z}^2} \in X_{A,B}^\kappa$
such that
$\omega_{0,0} =\omega,
 \omega_{i,j} =\omega'$
for some $(i,j) \in {\Bbb Z}^2$ with
$j<0<i$.
We set
\begin{equation}
\Omega_\kappa
=\{  (\alpha,a) \in   E_A\times E_B |  s(\alpha) = s(a),   
\kappa(\alpha,b) =(a,\beta) \text{ for some } \beta \in E_A, b\in E_B \} 
\label{eqn:OK}
\end{equation}
and define two $|\Omega_\kappa| \times |\Omega_\kappa|$-matrcies
$A_\kappa$ and $B_\kappa$  with entries in $\{0,1\}$
by 
\begin{align}
A_\kappa((\alpha,a),(\delta,b))
& = 
{\begin{cases}
1 &  \kappa(\alpha,b) = (a,\beta) \text{ for some } \beta\in E_A,\\  
0 & \text{ otherwise}
\end{cases}
} \label{eqn:Akappa}
\\
\intertext{ for
$(\alpha,a),(\delta,b) \in \Omega_\kappa$, and} 
B_\kappa((\alpha,a), (\beta,d))           
& = 
{\begin{cases}
1 & \kappa(\alpha,b) = (a,\beta) \text{ for some } b \in E_B,\\  
0 & \text{ otherwise}
\end{cases}
} \label{eqn:Bkappa}
 \end{align}
for 
$(\alpha,a), (\beta,d) \in \Omega_\kappa$
respectively.
Put the matrix
\begin{equation}
H_\kappa
=
\begin{bmatrix}
A_\kappa & A_\kappa \\
B_\kappa & B_\kappa 
\end{bmatrix}. \label{eqn:HK}
\end{equation}
It has been proved in
\cite{MaPre2011No2}
that the $C^*$-algebra ${\cal O}_{{\cal H}_\kappa^{A,B}}$
associated with the 
Hilbert $C^*$-quad module ${\cal H}_\kappa^{A,B}$
is isomorphic to the Cuntz-Krieger algebra
${\cal O}_{H_\kappa}$
for the matrix 
$
H_\kappa.
$
In this paper, we first show the following theorem.
\begin{thm}[Theorem \ref{thm:1.1}]
The subshift 
$X_{A,B}^\kappa$
of the tiling space 
is transitive if and if
the matrix $H_\kappa$ 
is irreducible.
In this case, $H_\kappa$ satisfies condition (I) in the sense of \cite{CK}. 
Hence 
if the subshift 
$X_{A,B}^\kappa$
of the tiling space 
is transitive,
 the $C^*$-algebra
$
{\cal O}_{{\cal H}_\kappa^{A,B}}
$
is simple and purely infinite.
\end{thm}
We will second see the following theorem.
 \begin{thm}[Theorem \ref{thm:1.2}]
If the matrix
 $A$ or $B$ is irreducible,
 the matrix
 $H_\kappa$ is irreducible and satisfies condition (I),
so that
 the $C^*$-algebra
$
{\cal O}_{{\cal H}_\kappa^{A,B}}
$
is simple and purely infinite.
\end{thm}
Let $N,M$ be positive integers with $N, M >1$.
They give 
$1 \times 1$ commuting matrices 
$A = [N], B=[M]$.
We will present  K-theory formulae 
for 
the $C^*$-algebras 
$
{\cal O}_{{\cal H}_\kappa^{[N],[M]}}
$
with exchanging specification
$\kappa$.
The directed graph $G_A$
associated to the matrix
$A =[N]$
is a graph consists of $N$-self directed loops denoted by $E_A$
with  a vertex denoted by $v$.
Similarly 
the directed graph $G_B$
consists of  $M$-self directed loops
denoted by $E_B$ with the vertex $v$.
We fix a specification
$\kappa: E_A \times E_B \longrightarrow E_B \times E_A$
defined by exchanging
$\kappa(\alpha,a) = (a,\alpha)$
for
$(\alpha,a) \in E_A \times E_B$.
We will have the following K-theory formulae for the $C^*$-algebra 
$
{\cal O}_{{\cal H}_\kappa^{[N],[M]}}.
$
In its computation, the Euclidean algorithm 
will be used.
\begin{thm}
For integers 
$1< N \le M \in{\Bbb N}$
and
a specification $\kappa$ 
of exchanging directed $N$-loops and $M$-loops,
the $C^*$-algebra
${\cal O}_{{\cal H}_\kappa^{[N],[M]}}
$
is a simple purely infinite Cuntz-Krieger algebra
whose K-groups are
\begin{align*}
K_1({\cal O}_{{\cal H}_\kappa^{[N],[M]}})
\cong 
& 0,\\
K_0({\cal O}_{{\cal H}_\kappa^{[N],[M]}}) 
\cong
&
\overbrace{
{\Bbb Z}/(N-1){\Bbb Z} \oplus \cdots \oplus {\Bbb Z}/(N-1){\Bbb Z}
}^{M-2} \\
\oplus
& \overbrace{
{\Bbb Z}/(M-1){\Bbb Z} \oplus \cdots \oplus {\Bbb Z}/(M-1){\Bbb Z}
}^{N-2} \\
\oplus
& 
{\Bbb Z}/ d {\Bbb Z}
\oplus
{\Bbb Z}/
[k_1,k_2,\dots,k_{j+1}] (M-1)(M+N-1)
{\Bbb Z}
 \end{align*}
where
$d = (N-1,M-1)$ the greatest common divisor of
$N-1$ and $M-1$,
and
the sequence
$k_0, k_2, \dots, k_{j+1}$
is the successive integral quotients of $M-1$ by $N-1$
by the Euclidean algorithm,   
and
the integer 
$[k_1,k_2,\dots, k_{j+1}]
$
is defined by inductively
\begin{align*}
[k_0]  = 1, \quad [k_1]  & = k_1, \quad
[k_1,k_2]  = 1 + k_1 k_2, \quad
              \dots,                   \\
[k_1,k_2,\dots, k_{j+1}] & = [k_1,k_2, \dots, k_{j}] k_{j+1} 
                           + [k_1,\dots,k_{j-1}]. 
\end{align*}
\end{thm}
We remark that the $C^*$-algebras studied in this paper
are different from the higher rank graph algebras
studied by 
A. Kumjian--D. Pask \cite{KP}, G. Robertson--T. Steger \cite{RoSte},
V. Deaconu \cite{Dea}, and etc. . 

Throughout the paper, 
we denote by ${\Bbb N}$ and by $\Zp$
the set of positive integers and the set of nonnegative integers respectively.


\section{Transitivity of the tilings $X_{A,B}^\kappa$ and 
simplicity of ${\cal O}_{{\cal H}_\kappa^{A,B}}$}


Let $\Sigma$ be a finite set.
 The two-dimensional full shift over $\Sigma$ is defined to be 
\begin{equation*}
\Sigma^{{\Bbb Z}^2}
= \{ (x_{i,j})_{(i,j) \in {\Bbb Z}^2} \mid x_{i,j} \in \Sigma \}.
\end{equation*}
An element $x \in \Sigma^{{\Bbb Z}^2}$
is regarded as a function
$x: {\Bbb Z}^2 \longrightarrow \Sigma
$
which is called a configuration on ${\Bbb Z}^2$.
For a vector $m=(m_1,m_2) \in {\Bbb Z}^2$,
let 
$
\sigma^m :  \Sigma^{{\Bbb Z}^2} \longrightarrow \Sigma^{{\Bbb Z}^2}
$
be
the translation along vector $m$ defined by
\begin{equation*}
\sigma^m ((x_{i,j})_{(i,j) \in {\Bbb Z}^2})
= (x_{i+{m_1}, j+{m_2}})_{(i,j) \in {\Bbb Z}^2}. 
\end{equation*}
A subset 
$X \subset \Sigma^{{\Bbb Z}^2}$ 
is said to be  translation invariant if
$\sigma^m(X) = X$ for all $m \in {\Bbb Z}^2$.
 It is obvious to see that 
a subset 
$X \subset \Sigma^{{\Bbb Z}^2}$ 
is  translation invariant if
 and only if 
 $X$ is invariant only both horizontaly and vertically,
 that is,
 $\sigma^{(1,0)}(X) = X$ and
$\sigma^{(0,1)}(X) = X$.
For $k \in \Zp$, 
put
 \begin{equation*}
[-k,k]^2 
= \{ (i,j) \in {\Bbb Z}^2 \mid -k\le i,j\le k \}
= [-k,k] \times [-k,k].
 \end{equation*}
A metric $d$ on
$\Sigma^{{\Bbb Z}^2}$ is defined by
for $x,y \in \Sigma^{{\Bbb Z}^2}$
with $x \ne y$ 
 \begin{equation*}
d(x,y)  
= 
\frac{1}{2^k}  \quad \text{ if } \quad x_{(0,0)} = y_{(0,0)},
 \end{equation*}
where
$k = \max \{k \in {\Bbb Z}_+ \mid  
x_{[-k,k]^2} = y_{[-k,k]^2} \}.
$ 
If $x_{(0,0)} \ne y_{(0,0)}$, put $k= -1$ on the above definition.
If $x = y$, we set $d(x,y) =0$.
A two-dimensional subshift $X$ is defined to be a closed,
translation invariant subset of $\Sigma^{{\Bbb Z}^2}$
(cf. \cite[p.467]{LM}).
A two-dimensional subshift $X$ is said to 
have the {\it diagonal property}\
if for 
$ (x_{i,j})_{(i,j)\in {\Bbb Z}^2}, (y_{i,j})_{(i,j)\in {\Bbb Z}^2}\in X$,
the conditions
$
x_{i,j} = y_{i,j},
x_{i+1,j-1} = y_{i+1,j-1}
$
imply 
$
x_{i,j-1} = y_{i,j-1},
x_{i+1,j} = y_{i+1,j}
$
(see \cite{MaPre2011}).
The diagonal property has the following property:
 for $x \in X$ and $(i,j) \in {\Bbb Z}^2$,
the configuration $x$ is determined by the diagonal line
$(x_{i+n, j-n})_{n \in {\Bbb Z}}$ 
through   
$(i,j)$.

We henceforth  go back to our previous situation of
$C^*$-textile dynamical system 
$(\A_N, \rho^A, \rho^B, E_A, E_B,\kappa)$
coming from $N \times N$ commuting matrices 
$A$ and $B$ with specification $\kappa$ as in Section 1.
We always assume that both matrices $A$ and $B$ are essential.
It yields 
a two-dimensional subshift 
$X_{A,B}^\kappa$
as follows:
Let $\Sigma $ be the set $E_\kappa$ 
of tiles defined in \eqref{eqn:EK}.
For 
$\omega= (\alpha, b, a,\beta) \in E_\kappa$,
define maps
$
t(=\text{top}), b(=\text{bottom}) :E_\kappa \longrightarrow E_A
$
and
$
l(=\text{left}), r(=\text{right}) :E_\kappa \longrightarrow E_B
$
by 
setting 
\begin{equation*}
t(\omega) = \alpha, \quad b(\omega) = \beta, \quad
l(\omega) = a, \quad r(\omega) = b
\end{equation*}
as in the following figure:
$$
\begin{CD}
\circ @>\alpha=t(\omega)>> \circ \\
@V{a=l(\omega)}VV  @VV{b=r(\omega)}V \\
\circ @>>\beta= b(\omega)> \circ 
\end{CD}
$$
A configuration
$(\omega_{i,j})_{(i,j)\in{\Bbb Z}^2} \in E_\kappa^{{\Bbb Z}^2}$
is said to be {\it paived}\
if the conditions 
\begin{equation*}
t(\omega_{i,j}) = b(\omega_{i,j+1}), \quad 
r(\omega_{i,j}) = l(\omega_{i+1,j}), \quad
l(\omega_{i,j}) = r(\omega_{i-1,j}), \quad
b(\omega_{i,j}) = t(\omega_{i,j-1})
\end{equation*}
hold for all $(i,j) \in {\Bbb Z}^2$.
Let
$X_{A,B}^\kappa$
be the set of all paved configurations
$(\omega_{i,j})_{(i,j) \in {\Bbb Z}^2} 
\in E_\kappa^{{\Bbb Z}^2}. 
$
It consists of the Wang tilings of the tiles of $E_\kappa$
(see \cite{Wang}).
The following proposition is easy.
\begin{prop}
$X_{A,B}^\kappa$ 
is a two-dimensional subshift having diagonal property.
\end{prop}
Let 
$e_\omega, \omega \in E_\kappa$ 
be the standard basis of ${\Bbb C}^{|E_\kappa|}$.
Put the projection
$E_\omega 
= \rho^B_b \circ \rho^A_\alpha(1)(= \rho^A_\beta \circ \rho^B_a(1))
\in \A_N$
for $\omega = (\alpha,b,a,\beta) \in E_\kappa$.
We set 
$$
{\cal H}_\kappa^{A,B} = \sum_{\omega \in E_\kappa} e_\omega \otimes E_\omega \A_N.
$$
Then 
$
{\cal H}_\kappa^{A,B}
$
has a natural structure of not only Hilbert $C^*$-right module over $\A_N$
but also two other Hilbert $C^*$-bimodule structure,
called Hilbert $C^*$-quad module.
By two-dimensional analogue of Pimsner's construction of Hilbert $C^*$-bimodule algebra
(\cite{Pim}),
we have introduced a $C^*$-algebra 
$
{\cal O}_{{\cal H}_\kappa^{A,B}}
$
(see \cite{MaPre2011No2} for detail construction).
Let $\Omega_\kappa$
be the subset of
$E_A \times E_B$ 
defined in
\eqref{eqn:OK}.
We define 
two $|\Omega_\kappa| \times |\Omega_\kappa|$-matrcies
$A_\kappa$ and $B_\kappa$  with entries in $\{0,1\}$
as in 
\eqref{eqn:Akappa} and  \eqref{eqn:Bkappa}.
The matrices 
$A_\kappa$ and $B_\kappa$
represent the concatenations of edges 
as in the following figures respectively:
\begin{align*}
{\begin{CD}
\circ @>\alpha>> \circ @>\delta>> \\
@V{a}VV  @V{b}VV @. \\ 
\circ @>>> \circ @.
\end{CD}
}
\qquad 
& \text{if } A_\kappa((\alpha,a),(\delta,b)) =1, \\
\intertext{ and }
 {
\begin{CD}
\circ @>\alpha>> \circ \\
@V{a}VV   @VVV \\
\circ @>\beta >> \circ  \\
@V{d}VV   @. 
\end{CD}
}
\qquad
\qquad
\qquad
& \text{if } B_\kappa((\alpha,a),(\beta,d)) =1.
\end{align*}
Let $H_\kappa$ be the $2|\Omega_\kappa| \times 2|\Omega_\kappa|$
matrix defined in \eqref{eqn:HK}.
We have proved 
the following result in
\cite{MaPre2011No2}.
\begin{thm}\label{thm:CK} 
The $C^*$-algebra ${\cal O}_{{\cal H}_\kappa^{A,B}}$
associated with 
Hilbert $C^*$-quad module ${\cal H}_\kappa^{A,B}$
defined by commuting matrices
$A, B$
and a specification $\kappa$
is isomorphic to
the Cuntz-Krieger algebra 
${\cal O}_{H_\kappa}$
for the matrix 
$H_\kappa$.
Its K-groups 
$K_*({\cal O}_{H_\kappa})$
are computed as 
\begin{align*}
K_0({\cal O}_{H_\kappa}) & 
={\Bbb Z}^{n} / (A_\kappa + B_\kappa - I_n){\Bbb Z}^n, \\
K_1({\cal O}_{H_\kappa}) & 
= \Ker(A_\kappa + B_\kappa - I_n) \text{ in } {\Bbb Z}^n,
\end{align*}
where
$n = |\Omega_\kappa|$.
\end{thm}
We will study 
a relationship between 
transitivity of the tiling space $X_{A,B}^\kappa$ and 
simplicity of the $C^*$-algebra ${\cal O}_{{\cal H}_\kappa^{A,B}}$.
An essential matrix with entries in $\{0,1\}$ 
is said to satisfy condition (I)
(in the sense of \cite{CK})
if the shift space defined by the topological Markov
chain for the matrix is homeomorphic to a Cantor discontinuum. 
The condition is equivalent to the condition that
every loop
in the associated directed graph 
has an exit (\cite{KPRR}).
It is a fundamental result that
a Cuntz-Krieger algebra is
simple and purely infinite if the underlying matrix is irreducible
and satisfies condition (I) (\cite{CK}). 
We will find a condition
of the two-dimensional subshift  
$X_{A,B}^\kappa$ 
of the tiling space
under which the matrix 
$H_\kappa$ is irreducible and satisfies condition (I).
Hence
the condition yields the simplicity and purely infiniteness 
of the algebra ${\cal O}_{{\cal H}_\kappa^{A,B}}$.

We are assuming that both of the matrices $A$ and $B$ are essential.
Then we have
\begin{lem}
Both of the matrices 
$A_\kappa$ and $B_\kappa$
are essential.
\end{lem}
\begin{pf}
For $(\alpha, a) \in \Omega_\kappa$,
by definition of $\Omega_\kappa$,
there exist $\beta \in E_A$ 
and
$b \in E_B$
such that
$\kappa(\alpha, b) = (a, \beta)$.
Since $A$ is essential,
one may take $\beta_1 \in E_A$
such that 
$s(\beta_1) = r(b)(=r(\beta)).
$
Hence
$(b, \beta_1) \in \Sigma^{BA}$.
Put
$
(\alpha_1,b_1) 
= \kappa^{-1}(b, \beta_1) \in \Sigma^{AB}
$
so that 
$(\alpha_1,b) \in \Omega_\kappa$
and
$A_\kappa((\alpha,a), (\alpha_1,b)) =1$
as in the following figure:
\begin{equation*}
\begin{CD}
\circ @>\alpha>> \circ @>\alpha_1>>\circ \\
@V{a}VV  @V{b}VV @V{b_1}VV \\
\circ @>\beta>> \circ @>\beta_1>>\circ 
\end{CD}
\end{equation*}
For
$(\delta,b) \in \Omega_\kappa$
there exists 
$\alpha \in E_A$ such that
$r(\alpha)= s(\delta)(=s(b))$
because 
$A$ is essential.
Hence
$(\alpha, b) \in \Sigma^{AB}$.
Put
$(a,\beta) = \kappa(\alpha, b)$
so that
$(\alpha,a) \in \Omega_\kappa$
and
$A_\kappa((\alpha,a),(\delta,b)) =1$
as in the following figure:
\begin{equation*}
\begin{CD}
\circ @>\alpha>> \circ @>\delta>> \\
@V{a}VV  @V{b}VV @.\\
\circ @>\beta>> \circ 
\end{CD}
\end{equation*}
Therefore one sees 
that $A_\kappa$ is essential,
and similarly 
that $B_\kappa$ is essential.
\end{pf}
Hence we have
\begin{prop}
The matrix $H_\kappa$ is essential and 
satisfies condition (I).
\end{prop}
\begin{pf}
By the previous lemma,
both of the matrices
$A_\kappa$ and $B_\kappa$ are essential.
Hence every row of $A_\kappa$ and of $B_\kappa$ 
has at least one $1$.
Since
$$
H_\kappa 
=
\begin{bmatrix}
A_\kappa & A_\kappa \\ 
B_\kappa & B_\kappa 
\end{bmatrix},
$$
every row of $H_\kappa$ has at least two $1'$s.
This implies that a loop in the directed graph associated to the matrix 
$H_\kappa$ must has an exit
so that
the matrix $H_\kappa$ 
satisfies condition (I).
\end{pf}

For 
$(\alpha,a), (\alpha',a') \in \Omega_\kappa$,
and 
$C, D =A$ or $B$,
we have
\begin{equation*}
[C_\kappa D_\kappa]((\alpha,a), (\alpha',a'))
= \sum_{(\alpha_1,a_1) \in \Omega_\kappa}
C_\kappa((\alpha,a), (\alpha_1,a_1))
D_\kappa((\alpha_1,a_1), (\alpha',a')).
\end{equation*}
Hence
$[A_\kappa A_\kappa]((\alpha,a), (\alpha',a')) \ne 0$
if and only if 
there exists
$(\alpha_1, a_1) \in \Omega_\kappa$
such that
$\kappa(\alpha,a_1) = (a,\beta) $ for some $\beta \in E_A$  
and
$\kappa(\alpha_1,a') = (a_1,\beta_1) $ for some $\beta_1 \in E_A$  
as in the following figure:
\begin{equation*}
\begin{CD}
\circ @>\alpha>> \circ @>\alpha_1>>\circ @>\alpha'>> \\
@V{a}VV  @V{a_1}VV  @V{a'}VV @.\\
\circ @>\beta>> \circ @>\beta_1>> \circ 
\end{CD}
\end{equation*}
And also
$[A_\kappa B_\kappa]((\alpha,a), (\alpha',a')) \ne 0$
if and only if 
there exists
$(\alpha_1, a_1) \in \Omega_\kappa$
such that
$\kappa(\alpha,a_1) = (a,\beta) $ for some $\beta \in E_A$  
and
$\kappa(\alpha_1,b_1) = (a_1,\alpha') $ for some $b_1 \in E_B$  
as in the following figure:
\begin{equation*}
\begin{CD}
\circ @>\alpha>> \circ @>\alpha_1>>\circ  \\
@V{a}VV  @V{a_1}VV  @V{b_1}VV @.\\
\circ @>\beta>> \circ @>\alpha'>> \circ \\
@.  @V{a'}VV @.
\end{CD}
\end{equation*}
Similarly 
$[B_\kappa A_\kappa]((\alpha,a), (\alpha',a')) \ne 0$
if and only if 
there exists
$(\alpha_1, a_1) \in \Omega_\kappa$
such that
$\kappa(\alpha,b) = (a,\alpha_1) $ for some $b \in E_B$  
and
$\kappa(\alpha_1,a') = (a_1,\beta_1) $ for some $\beta_1 \in E_A$  
as in the following figure:
\begin{equation*}
\begin{CD}
\circ @>\alpha>> \circ  \\
@V{a}VV  @V{b}VV   @.\\
\circ @>\alpha_1>> \circ @>\alpha'>> \\ 
@V{a_1}VV  @V{a'}VV   @.\\
\circ @>\beta_1>> \circ 
\end{CD}
\end{equation*}
And also
$[B_\kappa B_\kappa]((\alpha,a), (\alpha',a')) \ne 0$
if and only if 
there exists
$(\alpha_1, a_1) \in \Omega_\kappa$
such that
$\kappa(\alpha,b) = (a,\alpha_1) $ for some $b \in E_B$  
and
$\kappa(\alpha_1,b_1) = (a_1,\alpha') $ for some $b_1 \in E_B$  
as in the following figure:
\begin{equation*}
\begin{CD}
\circ @>\alpha>> \circ \\
@V{a}VV  @V{b}VV \\
\circ @>\alpha_1>> \circ \\
@V{a_1}VV  @V{b_1}VV \\
\circ @>\alpha'>> \circ \\
@V{a'}VV @.
\end{CD}
\end{equation*}
\begin{lem}
$A_\kappa B_\kappa = 
B_\kappa A_\kappa.
$
\end{lem}
\begin{pf}
For 
$(\alpha,a), (\alpha',a') \in \Omega_\kappa$,
we have
$[A_\kappa B_\kappa]((\alpha,a), (\alpha',a')) =m$
if and only if 
there exist
$(\alpha_i, a'_i) \in \Omega_\kappa, i=1,\dots,m$
such that
$\kappa(\alpha,a'_i) = (a,\beta_i) $ for some $\beta_i \in E_A$  
and
$\kappa(\alpha_i,b_i) = (a'_i,\alpha') $ for some $b_i \in E_B$  
as in the following figure:
\begin{equation*}
\begin{CD}
\circ @>\alpha>> \circ @>\alpha_i>>\circ  \\
@V{a}VV  @V{a'_i}VV  @V{b_i}VV @.\\
\circ @>\beta_i>> \circ @>\alpha'>> \circ \\
@.  @V{a'}VV @.
\end{CD}
\end{equation*}
Put
$(a_i, \beta'_i) = \kappa(\beta_i,a')$.
We then have
$(\beta_i,a_i) \in \Omega_\kappa $  
as in the following figure:
\begin{equation*}
\begin{CD}
\circ @>\alpha>> \circ  \\
@V{a}VV  @V{a'_i}VV   @.\\
\circ @>\beta_i>> \circ @>\alpha'>> \\ 
@V{a_i}VV  @V{a'}VV   @.\\
\circ @>\beta'_i>> \circ 
\end{CD}
\end{equation*}
If $(\beta_i,a_i) =(\beta_j,a_j)$ in $\Omega_\kappa$,
then we have $\beta_i = \beta_j$
so that $a'_i= a'_j$ and hence
$\alpha_i = \alpha_j$. 
Therefore we have 
$[B_\kappa A_\kappa]((\alpha,a), (\alpha',a')) =m$.
\end{pf}

\begin{lem}
The following four conditions are equivalent.
\begin{enumerate}
\renewcommand{\labelenumi}{(\roman{enumi})}
\item 
The matrix $H_\kappa$ is irreducible.
\item
For
$(\alpha,a), (\alpha',a') \in \Omega_\kappa$,
there exist
$n,m \in \Zp$ such that
\begin{equation*}
A_\kappa (A_\kappa + B_\kappa)^n((\alpha,a), (\alpha',a')) >0,
\qquad
B_\kappa (A_\kappa + B_\kappa)^m((\alpha,a), (\alpha',a')) >0.
\end{equation*}
\item
The matrix 
$A_\kappa + B_\kappa$ is irreducible.
\item
For
$(\alpha,a), (\alpha',a') \in \Omega_\kappa$,
there exists a paved configuration
$(\omega_{i,j})_{(i,j) \in {\Bbb Z}^2} \in X_{A,B}^\kappa$
such that
$$
t(\omega_{0,0}) =\alpha,\quad 
 l(\omega_{0,0}) =a, \quad
t(\omega_{i,j}) =\alpha', \quad
l(\omega_{i,j}) =a'
$$
for some $(i,j) \in {\Bbb Z}^2$ 
with 
$ j<0<i$.
\end{enumerate}
\end{lem}
\begin{pf}
(i) $\Longleftrightarrow$ (ii): 
The identity
\begin{equation}
H_\kappa^n =
\begin{bmatrix}
A_\kappa (A_\kappa + B_\kappa)^n & A_\kappa (A_\kappa + B_\kappa)^n \\
B_\kappa (A_\kappa + B_\kappa)^n & B_\kappa (A_\kappa + B_\kappa)^n
\end{bmatrix}
\end{equation}
implies the equivalence between 
(i) and (ii).

(ii) $\Longrightarrow$ (iii):
Suppose that for
$(\alpha,a), (\alpha',a') \in \Omega_\kappa$,
there exists $n \in \Zp$ such that
$
A_\kappa (A_\kappa + B_\kappa)^n ((\alpha,a), (\alpha',a')) >0$
so that
$$
  (A_\kappa + B_\kappa)^{n+1}((\alpha,a), (\alpha',a')) >0.
  $$
Hence
the matrix
$A_\kappa + B_\kappa$
is irreducible.

(iii) $\Longrightarrow$ (ii):
As
$A_\kappa$ and $B_\kappa$ are both essential,
for
$(\alpha, a), (\alpha',a') \in \Omega_\kappa$
there exists
$(\alpha_1, a_1), (\alpha_2,a_2) \in \Omega_\kappa$
such that
$$
A_\kappa((\alpha, a), (\alpha_1,a_1)) =1,
\qquad
B_\kappa((\alpha, a), (\alpha_2,a_2)) =1.
$$
Since $A_\kappa + B_\kappa$ is irreducible,
there exist
$n,m \in \Zp$ such that
\begin{equation*}
 (A_\kappa + B_\kappa)^n((\alpha_1,a_1), (\alpha',a')) >0,
\qquad
(A_\kappa + B_\kappa)^m((\alpha_2,a_2), (\alpha',a')) >0.
\end{equation*}
Hence
we have
\begin{equation*}
A_\kappa (A_\kappa + B_\kappa)^n((\alpha,a), (\alpha',a')) >0,
\qquad
B_\kappa (A_\kappa + B_\kappa)^m((\alpha,a), (\alpha',a')) >0.
\end{equation*}

(ii) $\Longrightarrow$ (iv):
For
$(\alpha, a), (\alpha',a') \in \Omega_\kappa$,
take
$(\alpha_1,a_1) \in \Omega_\kappa$ and
$\beta \in E_A$ 
such that 
$\kappa(\alpha,a_1) = (a,\beta)$.
By (ii), 
there exists $m \in \Zp$ with
$B_\kappa (A_\kappa + B_\kappa)^m((\alpha, a), (\alpha',a')) >0$.
One may take 
$b' \in E_B$
and 
$\beta' \in E_A$
 satisfying
$\kappa(\alpha',b') = (a',\beta')$,
so that
there exists a paved configuration 
$(\omega_{i,j})_{(i,j) \in {\Bbb Z}^2} \in X_{A,B}^\kappa$
such that
$\omega_{0,0} = (\alpha, a_1, a, \beta)$
and
$\omega_{i,j} = (\alpha', b', a', \beta')$
for some $(i,j) \in {\Bbb Z}^2$ with
$j<0<i$
as in the following figure:
\begin{equation*}
\begin{CD}
\circ @>\alpha>> \circ @>\alpha_1>>\circ   \\
@V{a}VV  @V{a_1}VV  @V{}VV @.  \\
\circ @>\beta>> \circ @>>> \ddots \\
@.  @V{}VV @.  
\end{CD}
\hspace{4cm}
\end{equation*}
\nopagebreak
\begin{equation*}
\hspace{5cm} \ddots 
\end{equation*}
\nopagebreak
\begin{equation*}
\hspace{6cm}
\begin{CD}
\circ @>\alpha'>> \circ    \\
@V{a'}VV  @V{b'}VV  \\
\circ @>\beta'>> \circ 
\end{CD}
\end{equation*}

(iv) $\Longrightarrow$ (ii):
  The assertion is clear.
\end{pf}

\noindent
{\bf Definition.}
A two-dimensional subshift
$X_{A,B}^\kappa$
is said to be {\it transitive}
if 
for
two tiles
$\omega, \omega' \in E_\kappa$
there exists a paved configuration
$(\omega_{i,j})_{(i,j) \in {\Bbb Z}^2} \in X_{A,B}^\kappa$
such that
$\omega_{0,0} =\omega
$ and 
$
 \omega_{i,j} =\omega'$
for some $(i,j) \in {\Bbb Z}^2$
with
$j<0<i$.

\begin{thm}
The subshift 
$X_{A,B}^\kappa$
of the tiling space is transitive if and only if
the matrix $H_\kappa$ 
is irreducible.
\end{thm}
\begin{pf}
Assume that the matrix $H_\kappa$ is irreducible.
Hence the condition (iv) in Lemma 2.6 holds.
Let
$
\omega = (\alpha, b, a, \beta),
\omega' = (\alpha', b', a', \beta') 
 \in E_\kappa$
be two tiles.
Since $A$ is essential,
there exists
$\beta_1 \in E_A$ 
such that
$r(\beta)( = r(b)) = s(\beta_1)$,
so that
$(b,\beta_1) \in \Sigma^{BA}$.
One may take
$(\alpha_1, b_1) \in \Sigma^{AB}$
such that
$\kappa(\alpha_1, b_1) = (b,\beta_1)$
and hence
$(\alpha_1,b) \in \Omega_\kappa$
as in the following figure:
\begin{equation*}
\begin{CD}
\circ @>\alpha>> \circ @>\alpha_1>>\circ \\
@V{a}VV  @V{b}VV @V{b_1}VV \\
\circ @>\beta>> \circ @>\beta_1>>\circ 
\end{CD}
\end{equation*}
For
$(\alpha_1, b), (\alpha', a') \in \Omega_\kappa$,
by (iv) in Lemma 2.6,
there exists
$(\omega_{i,j})_{(i,j) \in {\Bbb Z}^2} \in X_{A,B}^\kappa$
such that
$t(\omega_{0,0}) = \alpha_1, l(\omega_{0,0}) = b,
 t(\omega_{i,j}) = \alpha', l(\omega_{i,j}) = a'$
 for some
$(i,j) \in {\Bbb Z}^2$
with
$j<0<i$.
 Since
 $X_{A,B}^\kappa$ has diagonal property,
 there exists a paved configuration 
$(\omega'_{i,j})_{(i,j) \in {\Bbb Z}^2} \in X_{A,B}^\kappa$
such that
$\omega'_{0,0} = \omega,
\omega'_{i,j} = \omega'.
$
Hence
 $X_{A,B}^\kappa$ is transitive.

Conversely
assume that
 $X_{A,B}^\kappa$ is transitive.
  For $(\alpha,a), (\alpha',a')  \in \Omega_\kappa$,
there exist
$b,b' \in E_B$
and
$\beta, \beta' \in E_A$
such that
$\omega = (\alpha, b, a, \beta),
\omega' = (\alpha', b', a', \beta') \in E_\kappa$.
It is clear that the transitivity of $X_{A,B}^\kappa$
implies  the condition (iv) in Lemma 2.6,
so that 
$H_\kappa$ is irreducible.
\end{pf}

%

\begin{lem}
If $A$ or $B$ is irreducible,
$X_{A,B}^\kappa$ is transitive.
\end{lem}
\begin{pf}
Suppose that the matrix
$A$ is irreducible.
For two tiles
$\omega = (\alpha, b, a, \beta),
\omega' = (\alpha', b', a', \beta')
\in E_\kappa$,
there exist
 concatenated edges
$(\beta, \beta_1,\dots,\beta_n, \alpha')$
in the graph $G_A$ for some
edges
$\beta_1,\dots,\beta_n \in E_A$.
Since 
$X_{A,B}^\kappa$ has diagonal property,
there exists a configuration
$(\omega_{i,j})_{(i,j) \in {\Bbb Z}^2} \in X_{A,B}^\kappa$
such that
$\omega' = \omega_{i,j}$
for some
$i>0, j=-1$.
Hence
 $X_{A,B}^\kappa$ is transitive.
\end{pf}
Since
the $C^*$-algebra
${\cal O}_{{\cal H}_\kappa^{A,B}}
$
is isomorphic to the Cuntz-Krieger algebra
${\cal O}_{H_\kappa}$
by \cite{MaPre2011No2},
we see the following theorems.
\begin{thm}\label{thm:1.1}
The subshift 
$X_{A,B}^\kappa$
of the tiling space 
is transitive if and if
the matrix $H_\kappa$ 
is irreducible.
In this case, $H_\kappa$ satisfies condition (I). 
Hence 
if the subshift 
$X_{A,B}^\kappa$
of the tiling space 
is transitive,
 the $C^*$-algebra
$
{\cal O}_{{\cal H}_\kappa^{A,B}}
$
is simple and purely infinite.
\end{thm}
By Lemma 2.8, we have  
 \begin{thm}\label{thm:1.2}
If the matrix
 $A$ or $B$ is irreducible,
 the matrix
 $H_\kappa$ is irreducible and satisfies condition (I),
so that
 the $C^*$-algebra
$
{\cal O}_{{\cal H}_\kappa^{A,B}}
$
is simple and purely infinite.
\end{thm}



\section{ The algebra ${\cal O}_{{\cal H}_\kappa^{[N],[M]}}$
for two positive integers $N,M $}

Let $N,M$ be positive integers with $N, M >1$.
They give 
$1 \times 1$ commuting matrices 
$A = [N], B=[M]$.
We will present  K-theory formulae 
for 
the $C^*$-algebras 
$
{\cal O}_{{\cal H}_\kappa^{[N],[M]}}
$
with exchanging specification
$\kappa$.
In the computations below, we will use 
Euclidean algorithm to find order of the torsion part of the $K_0$-group. 
The directed graph $G_A$
for the matrix
$A =[N]$
is a graph consists of  $N$-self directed loops 
with a vertex denoted by $v$.
The   $N$-self directed loops are denoted by $E_A$.
Similarly 
the directed graph $G_B$ for $B=[M]$
consists of 
$M$-self directed loops
denoted by $E_B$
with
the vertex $v$.
We fix a specification
$\kappa: E_A \times E_B \longrightarrow E_B \times E_A$
defined by exchanging
$\kappa(\alpha,a) = (a,\alpha)$
for
$(\alpha,a) \in E_A \times E_B$.
Hence 
$\Omega_\kappa = E_A \times E_B$
so that 
$|\Omega_\kappa | = | E_A | \times | E_B| = N \times M$.
We then know
$A_\kappa((\alpha,a),(\delta,b)) = 1 $
if and only if 
$b=a$,
and
$B_\kappa((\alpha,a),(\beta,d)) = 1 $
if and only if 
$\beta=\alpha$ as in the following figures
respectively.
\begin{equation*}
{\begin{CD}
\circ @>\alpha>> \circ @>\delta>> \\
@V{a}VV  @V{a=b}VV @. 
\end{CD}
}
\qquad 
\text{and}
\qquad 
\begin{CD}
\circ @>\alpha>>  \\
@V{a}VV   @. \\
\circ @>\alpha=\beta >>  \\
@V{d}VV   @. 
\end{CD}
\end{equation*}
In \cite{MaPre2011No2},
the K-groups for the case
$N=2$ and $M=3$
have been computed such that
$$
K_0({\cal O}_{{\cal H}_\kappa^{[2],[3]}}) \cong {\Bbb Z}/8{\Bbb Z},
\qquad
K_1({\cal O}_{{\cal H}_\kappa^{[2],[3]}}) \cong 0.
$$
We will generalize the above computations.

Let
$I_n$ be the $n \times n$ identity matrix
and 
$E_n$ the $n \times n$ matrix whose entries are all $1'$s.
For an $N \times N$-matrix
$C = [c_{i,j}]_{i,j=1}^N$
and
an $M \times M$-matrix
$ D = [d_{k,l}]_{k,l=1}^M$,
denote by
$C \otimes D$
 the $NM \times NM$  matrix
 $$
C \otimes D
=
\begin{bmatrix}
c_{11}D & c_{12}D  & \dots & c_{1N}D  \\
c_{21}D & c_{22}D  & \dots & c_{2N}D  \\
\vdots  & \vdots   & \ddots& \vdots   \\
c_{N1}D & c_{N2}D  & \dots & c_{NN}D  
\end{bmatrix}.
$$
Hence we have
\begin{equation*}
E_N \otimes I_M =
\begin{bmatrix}
I_M    & I_M    & \dots & I_M    \\
I_M    & I_M    & \dots & I_M    \\
\vdots & \vdots & \ddots& \vdots \\
I_M    & I_M    & \dots & I_M 
\end{bmatrix},
\qquad
I_N \otimes E_M =
\begin{bmatrix}
E_M   & 0    & \dots & 0      \\
0     & E_M  & \ddots & \vdots \\
\vdots &\ddots      & \ddots& 0      \\
0     & \dots& 0     & E_M 
\end{bmatrix}.
\end{equation*}
We denote by
$E_{[N]}= \{\alpha_1,\dots, \alpha_N \}, 
 E_{[M]} =\{a_1,\dots, a_M \}$.
As 
$
\Omega_\kappa =E_{[N]} \times E_{[M]},
$
the basis of ${\Bbb C}^N \otimes {\Bbb C}^M$ 
are ordered  lexcographically  from left as in the following way:
\begin{equation}
 (\alpha_1,a_1),\dots, (\alpha_1, a_M),
   (\alpha_2,a_1),\dots, (\alpha_2, a_M),
\dots,
(\alpha_N,a_1),\dots, (\alpha_N, a_M) \label{eqn:basis} 
\end{equation}
Let
$A_\kappa$ and $B_\kappa$ be the matrices
defined in the previous section
for the matrices
$A = [N], B = [M]$
with exchanging specification $\kappa$.
The following lemma is direct.
\begin{lem}
The matrices 
$A_\kappa , B_\kappa$
are written as 
\begin{equation*}
A_\kappa = E_N \otimes I_M, \qquad
B_\kappa = I_N \otimes E_M
\end{equation*}
along the ordered basis \eqref{eqn:basis}.
Hence we have
\begin{equation}
A_\kappa + B_\kappa - I_{NM}
 =
\begin{bmatrix}
E_M     & I_M      & \dots   & I_M     \\
I_M      & \ddots & \ddots & \vdots \\
\vdots & \ddots & \ddots & I_M      \\
I_M     & \dots   & I_M      & E_M 
\end{bmatrix}. \label{eqn:ABI}
\end{equation}
\end{lem}
We denote by 
$H_0$ the matrix
$A_\kappa + B_\kappa - I_{NM}$.
By Theorem 2.2, the K-groups of the algebra
$
{\cal O}_{{\cal H}_\kappa^{[N],[M]}}
$
are given by the kernel
$\Ker(H_0)$ 
and the cokernel
$\Coker(H_0)$
of the matrix $H_0$ in
${\Bbb Z}^{NM}$.
We will transform $H_0$ preserving isomorphism
classes of the groups 
$\Ker(H_0)$ and 
$\Coker(H_0)$ in ${\Bbb Z}^{NM}$
by the following operations called 
elementary operations on the matrix.

(A) Exchange  two rows or two columns.

(B) Multiply a row or column by $-1$.

(C) Add an integer multiple of one row to another row,
  or  of one column to another column. 
  
(D) Add a row vector 
    obtained by multiplication of an invertible matrix over ${\Bbb Z}$ of one row 
  to another row, or  of one column to another column. 

The
isomorphism classes of the groups of its 
kernel and its cokernel
do not change 
under the
elementary operations on the matrix.
We will successively apply the above elementary operations to the matrix $H_0$ 
to obtain  a diagonal matrix as in the following way.

(1) Add the minus of the $(i+1)$-th row to the $i$-th row in order for $i=1,\dots, N-1$
in $H_0$ to obtain the matrix below denoted by $H_1:$   

\begin{equation*}
H_1 =
\begin{bmatrix}
E_M-I_M & I_M-E_M  & 0        & \dots     & 0 \\
0          & E_M-I_M  & \ddots & \ddots   & \vdots \\
\vdots   & \ddots   & \ddots &  I_M-E_M & 0  \\
0          & \dots     &  0       & E_M-I_M  & I_M-E_M \\
I_M       & \dots     & \dots   & I_M        & E_M
\end{bmatrix}.
\end{equation*}

(2) Add the $i$-th row to the $(i+1)$-th row in order for $i=1,\dots, N-1$
in $H_1$
to obtain the matrix below denoted by $H_2:$  

\begin{equation*}
H_2 =
\begin{bmatrix}
E_M-I_M & I_M-E_M & 0      & \dots   & 0 \\
\vdots   & 0       & \ddots & \ddots  & \vdots \\
\vdots  & \vdots  & \ddots & I_M-E_M & 0  \\
E_M-I_M & 0       & \dots  & 0       & I_M-E_M \\
E_M      & I_M     & \dots  & I_M     & I_M
\end{bmatrix}.
\end{equation*}

(3) Add the $E_M - I_M$ multiplication  of the $N$-th row to the $(N-1)$-th row
in $H_2$
to obtain the matrix below denoted by $H_3:$  

\begin{equation*}
H_3 =
\begin{bmatrix}
E_M-I_M    & I_M-E_M  & 0        & \dots    & \dots   & 0 \\
\vdots      & 0           & \ddots & \ddots  &            & \vdots \\
\vdots      & \vdots   & \ddots & \ddots   & \ddots  & \vdots  \\
E_M-I_M    & 0           & \dots  & 0          & I_M-E_M & 0 \\
E_M^2-I_M & E_M-I_M & \dots  & \dots     & E_M-I_M & 0 \\
E_M         & I_M        & \dots   & \dots    & I_M        & I_M
\end{bmatrix}.
\end{equation*}

(4) Add the $i$-th row to the $(i-1)$-th row in order for $i=N-1,\dots, 2$
in $H_3$
to obtain  the matrix below denoted by $H_4:$  
\begin{equation*}
H_4 =
\begin{bmatrix}
p_M(N-1)       & 0          &          & \dots    &            & 0         \\
p_M(N-2)       & E_M-I_M &         &             &            &            \\
\vdots         & \vdots  & \ddots & \ddots   &            & \vdots  \\
p_M(2)         & \vdots  &           & \ddots   &            &            \\
p_M(1)         & E_M-I_M &\dots   & \dots   & E_M-I_M & 0        \\
E_M            & I_M        & \dots    & \dots    & I_M     & I_M
\end{bmatrix}
\end{equation*}
where
$p_M(i) = E_M^2+ (i-1) E_M- i I_M =(E_M + iI_M)(E_M - I_M)$ 
for $i=1,\dots N-1$.

(5) Add the minus of the $j$-th column 
to the  $(j-1)$-th column in order for $j=N,\dots, 3$
and 
the $-E_M$ multiplication of 
 the second column  to the first column 
in $H_4$
to obtain the matrix below denoted by $H_5:$  
\begin{equation*}
H_5 =
\begin{bmatrix}
p_M(N-1)      & 0          &          & \dots      &            & 0 \\
p_M(N-2)      & E_M-I_M &          &              &             &        \\
\vdots         & \vdots  & \ddots & \ddots    &             & \vdots  \\
p_M(2)         & \vdots  &           & \ddots    &             &    \\
p_M(1)         & E_M-I_M &\dots   & \dots     & E_M-I_M & 0 \\
0                &   0        & \dots   & \dots     & 0         & I_M
\end{bmatrix}.
\end{equation*}

(6)
 Add the minus of the $(N-1)$-th column 
to the $j$-th column in order for $j=N-2,\dots, 2$
and 
the $-(E_M + I_M)$ multiplication of the $(N-1)$-th column 
to the first column 
in $H_5$
to obtain the matrix below denoted by $H_6:$  
\begin{equation*}
H_6 =
\begin{bmatrix}
p_M(N-1)   & 0          &          & \dots    &         & 0 \\
p_M(N-2)   & E_M-I_M &          &            &         &        \\
\vdots       & 0         & \ddots & \ddots  &         & \vdots  \\
p_M(2)       & \vdots  & \ddots & \ddots  &         &   \\
0              & \vdots   &           & \ddots  & E_M-I_M & 0 \\
0              & 0         & \dots  & \dots    & 0       & I_M
\end{bmatrix}.
\end{equation*}

(7)
 Add
the  $-(E_M + (N-j) I_M)$ multiplication of the $j$-th column 
to the first column in order for $j= N-1,\dots, 2$
in $H_6$
to obtain the diagonal matrix below denoted by $H_7:$  
\begin{equation*}
H_7 =
\begin{bmatrix}
p_M(N-1)   & 0           &           & \dots    &            & 0 \\
0              & E_M-I_M &           &             &            &        \\
                &            & \ddots & \ddots   &            & \vdots  \\
\vdots       &            & \ddots & \ddots    &            &        \\
                &            &           &              & E_M-I_M & 0 \\
0              &            & \dots   &              & 0          & I_M
\end{bmatrix}.
\end{equation*}
As
$E_M^2 = M E_M$,
we have
$p_M(N-1) = (M+N-2)E_M -(N-1)I_M$.
We thus have
\begin{lem}
\begin{align*}
& \Ker (A_\kappa + B_\kappa - I_{NM}) \text{ in } {\Bbb Z}^{NM}
\cong 0 \\
\intertext{ and }
& \Coker (A_\kappa + B_\kappa - I_{NM}) \text{ in } {\Bbb Z}^{NM} \\
\cong
&
\overbrace{
{\Bbb Z}^M / (E_M -I_M){\Bbb Z}^M \oplus
\cdots \oplus {\Bbb Z}^M / (E_M -I_M){\Bbb Z}^M
}^{(N-2)}\\
\oplus 
&
{\Bbb Z}^M / ((M+N-2)E_M -(N-1)I_M){\Bbb Z}^M. 
\end{align*}
\end{lem}
\begin{pf}
It is straightforward to see that the matrix
$A_\kappa + B_\kappa - I_{NM}$
is invertivle by the formula \eqref{eqn:ABI}.
Since
$$
\Coker (A_\kappa + B_\kappa - I_{NM}) \text{ in } {\Bbb Z}^{NM}
\cong 
{\Bbb Z}^{NM} / H_7 {\Bbb Z}^{NM},
$$
the formula for the cokernel
is obvious.
\end{pf}
We will next compute the following groups
to compute 
$\Coker (A_\kappa + B_\kappa - I_{NM}) \text{ in } {\Bbb Z}^{NM}$.

(i) $ {\Bbb Z}^M / (E_M -I_M){\Bbb Z}^M,$ 

(ii) ${\Bbb Z}^M / ((M+N-2)E_M -(N-1)I_M){\Bbb Z}^M $

\medskip

(i)
As the matrix
$E_M -I_M$
is of the form
\begin{equation*}
\begin{bmatrix}
0        & 1         & \dots  & 1        \\
1        & 0         & \ddots & \vdots   \\
\vdots & \ddots & \ddots & 1 \\
1        & \dots & 1      & 0 
\end{bmatrix}.
\end{equation*}
by the same operations (1), (2) to get the matrix 
$H_2$ from $H_0$,
the matrix $E_M -I_M$
goes to the matrix
\begin{equation*}
\begin{bmatrix}
 -1    & 1      & 0      & \dots   & 0 \\ 
 -1    & 0      & \ddots & \ddots  & \vdots \\
\vdots & \vdots & \ddots & 1       & 0  \\
 -1    & 0      & \dots  & 0       & 1   \\
0      & 1      & \dots  & 1       & 1
\end{bmatrix}.
\end{equation*}
Add the minus of the $i$-th row to the $M$-th row 
in order for
$i =1,\dots,M-1$, 
we have the matrix
\begin{equation*}
\begin{bmatrix}
 -1    & 1      & 0      & \dots   & 0 \\ 
 -1    & 0      & \ddots & \ddots  & \vdots \\
\vdots & \vdots & \ddots & 1       & 0  \\
 -1    & 0      & \dots  & 0       & 1   \\
M-1   & 0      & \dots  & 0       & 0
\end{bmatrix}.
\end{equation*}
Add the $j$-th column to the first column
 for $j=2,\dots,M$, 
we have the matrix
\begin{equation*}
\begin{bmatrix}
 0     & 1      & 0      & \dots   & 0 \\ 
 0     & 0      & \ddots & \ddots  & \vdots \\
\vdots & \vdots & \ddots & 1       & 0  \\
 0     & 0      & \dots  & 0       & 1   \\
M-1   & 0      & \dots  & 0       & 0
\end{bmatrix} 
\end{equation*} 
which 
goes to 
the diagonal matrix with diagonal entries
$[1,1,\dots, 1, M-1]$
by exchanging rows.
Hence we see that
\begin{equation} 
{\Bbb Z}^M / (E_M -I_M){\Bbb Z}^M 
\cong 
{\Bbb Z} / (M-1) {\Bbb Z}. \label{eqn:(i)}
\end{equation}

(ii)
Put
$e = (M+N-2) -(N-1) = M-1$
and
$f = M + N -2$.
Then we have
\begin{equation}
(M+N-2)E_M -(N-1)I_M
=
\begin{bmatrix}
e        & f        & \dots  & f        \\
f         & e       & \ddots & \vdots   \\
\vdots & \ddots & \ddots & f \\
f         & \dots & f      & e 
\end{bmatrix}. \label{eqn:(ii)}
\end{equation}
By a similar manner  to the preceding operations from 
$H_1$ to $H_5$,
one obtains the following matrix denoted by $L_2$ 
from the matrix \eqref{eqn:(ii)}
\begin{equation*}
L_2 =
\begin{bmatrix}
e-f    & f-e    & 0      & \dots   & 0 \\
e-f    & 0      & \ddots & \ddots  & \vdots \\
\vdots & \vdots & \ddots & f-e     & 0  \\
e-f    & 0      & \dots  & 0       & f-e \\
e      & f      & \dots  & f       & f
\end{bmatrix}.
\end{equation*}
Add the $j$-th column to the first column for 
$j=2,\dots,M$
to obtain the matrix below denoted by $L_3:$
\begin{equation*}
L_3 =
\begin{bmatrix}
0        & f-e    & 0      & \dots   & 0 \\
0        & 0      & \ddots & \ddots  & \vdots \\
\vdots   & \vdots & \ddots & f-e     & 0  \\
0        & 0      & \dots  & 0       & f-e \\
e+(M-1)f & f      & \dots  & f       & f
\end{bmatrix}.
\end{equation*}
Exchang columns 
to obtain the matrix below denoted by $L_4:$
\begin{equation*}
L_4 =
\begin{bmatrix}
f-e    & 0      & \dots   & 0      & 0 \\
0      & \ddots & \ddots  & \vdots & \vdots \\
\vdots & \ddots & f-e     & 0      & 0 \\
0      & \dots  & 0       & f-e    & 0 \\
f      & \dots  & f       & f      &e+(M-1)f 
\end{bmatrix}.
\end{equation*}
Add the minus of the $j$-th column to the $(j-1)$-th column in order for 
$j=2,\dots,M-1$
to obtain the matrix below denoted by $L_5:$
\begin{equation*}
L_5 =
\begin{bmatrix}
f-e      & 0         &   \dots    & \dots       & 0 \\
e-f      & \ddots  & \ddots    &               & \vdots \\
0        & \ddots  & \ddots    & \ddots      & \vdots \\
\vdots & \ddots  & e-f         & f-e          & 0 \\
0        & \dots   & 0            & f              &e+(M-1)f 
\end{bmatrix}.
\end{equation*}
Add the $i$-th row to the $i+1$-th row in order for 
$i=1,\dots,M-2$
to obtain the matrix below denoted by $L_6:$
\begin{equation*}
L_6 =
\begin{bmatrix}
f-e      & 0          & \dots   & \dots   & 0 \\
0        & \ddots   & \ddots  &           & \vdots \\
\vdots & \ddots  & \ddots   & \ddots  & \vdots \\
0        & \dots    & 0        & f-e       & 0 \\
0        & \dots    & 0        & f          &e+(M-1)f 
\end{bmatrix}.
\end{equation*}
Put the $2 \times 2$ matrix $L_{(N,M)}$ 
by setting
\begin{equation*}
L_{(N,M)} =
\begin{bmatrix}
f-e    & 0 \\
f      &e+(M-1)f 
\end{bmatrix}.
\end{equation*}
As
$f -e =N-1$,
we have the following lemma with \eqref{eqn:(i)}.
\begin{lem} \hspace{6cm}
\begin{enumerate}
\renewcommand{\labelenumi}{(\roman{enumi})}
\item 
$ \qquad {\Bbb Z}^M / (E_M -I_M){\Bbb Z}^M 
\cong
{\Bbb Z} / (M-1) {\Bbb Z}.$ 
\item
 $\qquad {\Bbb Z}^M / ((M+N-2)E_M -(N-1)I_M){\Bbb Z}^M$ 
 
$ \cong
\overbrace{{\Bbb Z} / (N-1) {\Bbb Z}
\oplus \cdots \oplus {\Bbb Z} / (N-1) {\Bbb Z}
}^{M-2}
\oplus
{\Bbb Z}^2 / L_{(N,M)}{\Bbb Z}^2.$ 
 \end{enumerate}
 \end{lem}
It remains to compute the group
${\Bbb Z}^2 / L_{(N,M)}{\Bbb Z}^2.$ 
Put
$n = N-1, m= M-1$.
As
$ f - e = n$ and 
$f = m+n $,
we have
$e + (M-1) f = (M-1)(M + N-1) = m(m+n+1)$
so that
$$
L_{(N,M)} =
\begin{bmatrix}
n    & 0 \\
n+m  & m(m+n+1) 
\end{bmatrix}.
$$
Add the minus of the first row to the second row
in $L_{(N,M)}$ to obtain the matrix below denoted by $L_{n,m}:$ 
\begin{equation*}
L_{n,m} =
\begin{bmatrix}
n    & 0 \\
m  & m(m+n+1) 
\end{bmatrix}.
\end{equation*}
We may assume that
$M  \ge N$ and hence $m \ge n$.

If $m$ is divided by $n$  and hence there exists $k \in {\Bbb N}$ such that 
$m = nk$,
by adding  the  $-k$ multiplication of the first row to the second row
in $L_{n,m}$,
the matrix goes to the diagonal matrix: 
\begin{equation*}
\begin{bmatrix}
n    & 0 \\
0    & m(m+ n+1) 
\end{bmatrix}
=
\begin{bmatrix}
N-1    & 0 \\
0      & (M-1)(M + N -1) 
\end{bmatrix}.
\end{equation*}
Hence we have
\begin{equation*}
{\Bbb Z}^2 / L_{(N,M)}{\Bbb Z}^2
\cong
{\Bbb Z}/(N-1){\Bbb Z}
\oplus
{\Bbb Z}/(M-1)(M + N -1){\Bbb Z}.
\end{equation*}
Otherwise,
 by the Euclidean algorithm, we have 
 lists  of integers
 $r_0, r_1, \dots, r_j$
 and
 $k_0, k_1, \dots, k_{j+1}$
 for some $j \in {\Bbb N}$
 such that
\begin{align*}
m & = n k_0 + r_0, \qquad 0 < r_0 <n,   \\
n & = r_0 k_1 + r_1, \qquad 0 < r_1 <r_0,  \\
r_0 & = r_1 k_2 + r_2, \qquad 0 < r_2 <r_1,  \\
    & \dots                                 \\
r_{j-2} & = r_{j-1} k_j + r_j, \qquad 0 < r_j <r_{j-1},  \\
r_{j-1} & = r_j k_{j+1}, \qquad \quad 0 = r_{j+1}
\end{align*} 
where
$r_j = (m,n)$ the greatest common divisor of $m$ and $n$.
Put
$g = m(m+n+1)$.
Add the $-k_0$ multiplication of the first row to the second row
in $L_{n,m}$
to obtain the matrix below denoted by $L_{n,m}(0)$:
\begin{equation*}
L_{n,m}(0)
=
\begin{bmatrix}
n      & 0 \\
r_0    & g 
\end{bmatrix}.
\end{equation*}
Add the $-k_1$ multiplication of the second row to the first row
in $L_{n,m}(0)$
to obtain the matrix below denoted by $L_{n,m}(1)$:
\begin{equation*}
L_{n,m}(1)
=\begin{bmatrix}
r_1    & - k_1 g \\
r_0    & g 
\end{bmatrix}.
\end{equation*}
Add the $-k_2$ multiplication of the first row to the second row
in $L_{n,m}(1)$
to obtain the matrix below denoted by $L_{n,m}(2)$:
\begin{equation*}
L_{n,m}(2)
=
\begin{bmatrix}
r_1    & -k_1 g \\
r_2    & (1 + k_1 k_2)g 
\end{bmatrix}.
\end{equation*}
We continue these procedures as follows.
Add the $-k_{2i-1}$ multiplication of the second row to the first row
in $L_{n,m}(2i-2)$
to obtain the matrix dnoted by $L_{n,m}(2i-1)$.
And 
add the $-k_{2i}$ multiplication of  the first row to the second row
in $L_{n,m}(2i-1)$
to obtain the matrix denoted by $L_{n,m}(2i)$
for
$i=1,2,\dots$.
The algorithm stops at $j+1=2i-1$ or $j+1= 2 i$ for some $i \in {\Bbb N}$.
We set
\begin{align*}
[k_0]  = 1,\quad
[k_1]  & = k_1, \quad
[k_1,k_2]  = 1 + k_1 k_2, \quad
[k_1,k_2,k_3]  = [k_1,k_2] k_3 + [k_1],\quad \dots \\             
[k_1,k_2,\dots, k_{j+1}] & = [k_1,k_2, \dots, k_{j}] k_{j+1} + [k_1,\dots,k_{j-1}]. 
\end{align*}
Then we have
\begin{equation*}
L_{n,m}(1)
=\begin{bmatrix}
r_1    & - [k_1] g \\
r_0    & g 
\end{bmatrix},
\qquad
L_{n,m}(2)
=
\begin{bmatrix}
r_1    & -[k_1] g \\
r_2    & [k_1, k_2] g 
\end{bmatrix},
\end{equation*}
and inductively 
\begin{align*}
L_{n,m}(2i-1)
& =
\begin{bmatrix}
r_{2i-1}    & -[k_1,k_2,\dots,k_{2i-1}] g \\
r_{2i-2}    & [ k_1, k_2,\dots,k_{2i-2}] g 
\end{bmatrix},\\
L_{n,m}(2i)
& =
\begin{bmatrix}
r_{2i-1}    & -[k_1, k_2, \dots,k_{2i-1}] g \\
r_{2i}      & [k_1, k_2,\dots,k_{2i}] g 
\end{bmatrix}
\end{align*}
for $i=1,2,\dots $.
We denote by 
$d$ the greatest common divisor $(m,n)$
of $m$ and $n$,
so that $d = r_j$.
Take $m_0 \in {\Bbb Z}$ 
such that
$ m = m_0 d. $
Put
$
g_0 = m_0(m + n + 1) 
$
so that
$
 g = g_0 d.
$

We have two cases.

Case 1: $j+1 = 2i -1$ for some $i \in {\Bbb N}$. 
We have
\begin{equation*}
L_{n,m}(j+1)
=
\begin{bmatrix}
r_{j+1}    & -[k_1,k_2,\dots,k_{j+1} ] g \\
r_{j}      & [ k_1, k_2,\dots,k_j ] g 
\end{bmatrix}
=
\begin{bmatrix}
0    & -[k_1,k_2,\dots,k_{j+1}] g \\
d    & [ k_1, k_2,\dots,k_j ] g_0 d 
\end{bmatrix}.
\end{equation*}
Add the  $-[ k_1, k_2,\dots,k_{j}] g_0$ 
multiplication of the first column to the second column
in the above matrix $L_{n,m}(j+1)$,
and then exchange the rows
to obtain the matrix below
\begin{equation*}
\begin{bmatrix}
d    & 0 \\
0    & -[k_1,k_2,\dots,k_{j+1}] g
\end{bmatrix}.
\end{equation*}

Case 2: $j+1 = 2i$ for some $i \in {\Bbb N}$.
We have
\begin{equation*}
L_{n,m}(j+1)
=
\begin{bmatrix}
r_{j}    & -[k_1,k_2,\dots,k_{j}] g \\
r_{j+1}  & [ k_1, k_2,\dots,k_{j+1}] g 
\end{bmatrix}
=
\begin{bmatrix}
d    & - [ k_1, k_2,\dots,k_{j}] g_0 d \\ 
0    &   [k_1,k_2,\dots,k_{j+1}] g \\
\end{bmatrix}.
\end{equation*}
Add the $[ k_1, k_2,\dots,k_{j}] g_0$ 
multiplication of the first column to the second column
in the above matrix $L_{n,m}(j+1)$,
and then exchange the rows
to obtain the matrix below
\begin{equation*}
\begin{bmatrix}
d    & 0 \\
0    & [k_1,k_2,\dots,k_{j+1}] g
\end{bmatrix}.
\end{equation*}
We reach  the following lemma.
\begin{lem}
\begin{equation*}
{\Bbb Z}^2 / L_{(N,M)}
{\Bbb Z}^2
\cong 
{\Bbb Z}/ d {\Bbb Z}
\oplus
{\Bbb Z}/
[k_1,k_2,\dots,k_{j+1}] g
{\Bbb Z}.
\end{equation*}
\end{lem}
Therefore we have
\begin{thm}
For positive integers 
$1< N \le M \in{\Bbb N}$
and
a specification $\kappa$ 
of exchanging  $N$-loops and $M$-loops
in a  graph with one vertex,
the $C^*$-algebra
${\cal O}_{{\cal H}_\kappa^{[N],[M]}}
$
is a simple purely infinite Cuntz-Krieger algebra
whose K-groups are
\begin{align*}
K_1({\cal O}_{{\cal H}_\kappa^{[N],[M]}})\cong 
& 0,\\
K_0({\cal O}_{{\cal H}_\kappa^{[N],[M]}}) \cong
&
\overbrace{
{\Bbb Z}/(N-1){\Bbb Z} \oplus \cdots \oplus {\Bbb Z}/(N-1){\Bbb Z}
}^{M-2} \\
\oplus 
& 
\overbrace{
{\Bbb Z}/(M-1){\Bbb Z} \oplus \cdots \oplus {\Bbb Z}/(M-1){\Bbb Z}
}^{N-2} \\
\oplus
& 
{\Bbb Z}/ d {\Bbb Z}
\oplus
{\Bbb Z}/
[k_1,k_2,\dots,k_{j+1}] (M-1)(M+N-1)
{\Bbb Z}
 \end{align*}
where
$d = (N-1,M-1)$ is the greatest common divisor of
$N-1$ and $M-1$,
the sequence
$k_0, k_2, \dots, k_{j+1}$
of integers
is the list of the successive integral quotients of $M-1$ by $N-1$
in the Euclidean algorithm,   
and
the integer 
$[k_1,k_2,\dots, k_{j+1}]
$
is defined by inductively
\begin{align*}
[k_0]  = 1,\qquad [k_1] &  = k_1, \qquad
[k_1,k_2]  = 1 + k_1 k_2, \qquad \dots        \\
[k_1,k_2,\dots, k_{j+1}] & = [k_1,k_2, \dots, k_{j}] k_{j+1} 
                           + [k_1,\dots,k_{j-1}]. 
\end{align*}
\end{thm}
For the case 
$N=2 $ and $M\ge 2$,
we have
$d = 1,  r_0 =0$.
We understand 
$[k_1,\dots,k_{j+1}] =[k_0] = 1$
so that
we have
$$ 
[k_1,\dots,k_{j+1}] (M-1)(M + N -1) = 1 \times (M-1)(M+1) = M^2 -1.
$$
Hence
\begin{equation}
K_0({\cal O}_{{\cal H}_\kappa^{[2],[M]}}) \cong
{\Bbb Z}/ (M^2 -1) {\Bbb Z}. \label{eqn:35}
\end{equation}
If in particular, $M=3$,
the formula \eqref{eqn:35}
 is already seen in \cite{MaPre2011No2}.


\begin{thebibliography}{99}


%











%






 










\bibitem{C}
{\sc J. Cuntz},
{\it Simple $C^*$-algebras generated by isometries},
 Comm.\ Math.\ Phys.\
{\bf 57}(1977) pp.\ 173--185.










\bibitem{CK}{\sc J. ~Cuntz and W. ~Krieger},
{\it A class of $C^*$-algebras and topological Markov chains},
 Invent.\ Math.\
 {\bf 56}(1980), pp.\ 251--268.

\bibitem{Dea}{\sc V. ~Deaconu},
{\it  $C^*$-algebras and Fell bundles associated to a textile system},
preprint, arXiv:1001.0037.
























\bibitem{EGS}
{\sc R. ~Excel, D. ~Gon\c{c}alves and C.~Starling},
{\it The tiling $C^*$-algebra viewed as a tight inverse semigroup algebra},
 preprint arXiv:1106.4535v1.



\bibitem{EL}
{\sc R. Excel and M. Laca},
{\it Cuntz--Krieger algebras for infinite matrices},
J. Reine Angew. Math.
{\bf 512}(1999), pp.\ 119--172.
























\bibitem{KPW}
{\sc T. Kajiwara, C. Pinzari and Y. Watatani}, 
{\it Ideal structure and simplicity of the $C^*$-algebras generated by Hilbert modules},
J. Funct. Anal. {\bf 159}(1998) pp. 295--322.



































 
 


\bibitem{KP}
{\sc A. Kumjian and D. Pask},
{\it  Higher rank graph  $C^*$-algebras},
New York J.\ Math.\
{\bf 6}(2000), pp.\ 1-20.





\bibitem{KPRR}{\sc A. ~Kumjian, D. ~Pask, I. ~Raeburn and J. ~Renault}, 
{\it Graphs, groupoids and Cuntz--Krieger algebras},
 J. Funct. Anal.
{\bf 144}
(1997)
pp.\ 505--541.









\bibitem{LM}{\sc D. ~Lind and B. ~Marcus},
{\it An introduction to symbolic dynamics and coding},
 Cambridge University Press, Cambridge UK
(1995).



\bibitem{MP}
{\sc N. G.  ~Markley and M. E. ~Paul},
{\it Matrix subshifts for ${\Bbb Z}^\nu$ symbolic dynamics},
Proc.\ London  Math.\
{\bf 43}(1981), pp.\  251--272.







  
  





















\bibitem{MaCrelle}
{\sc K. Matsumoto},
{\it  Actions of symbolic dynamical systems on $C^*$-algebras},
J.\ Reine Angew.\ Math.\
{\bf 605}(2007), pp.\ 23-49.














\bibitem{MaMZ2010}
{\sc K. Matsumoto},
{\it  Actions of symbolic dynamical systems on $C^*$-algebras II.
Simplicity of $C^*$-symbolic crossed products and some examples},
Math.\ Z. 
{\bf 265}(2010), pp.\ 735-760.


\bibitem{MaPre2011}
{\sc K. Matsumoto},
{\it  $C^*$-algebras associated with textile dynamical systems},
preprint, arXiv:1106.5092v1. 

\bibitem{MaPre2011No2}
{\sc K. Matsumoto},
{\it  $C^*$-algebras associated with Hilbert $C^*$-quad modules of $C^*$-textile dynamical systems},
preprint, arXiv:1111.3091v1. 


\bibitem{Ma2011}
{\sc K. Matsumoto},
{\it  $C^*$-algebras associated with Hilbert $C^*$-multi modules (tentative title), 
in preparation}.







\bibitem{NaMemoir}
{\sc M. Nasu},
{\it Textile systems for endomorphisms and automorphisms of the shift}, 
Mem.\ Amer.\ Math.\ Soc.\ {\bf 546}(1995).







\bibitem{PRW1}
{\sc D. ~Pask, I.~Raeburn and N. A. ~Weaver},
{\it A family of $2$-graphs arising from two-dimensional subshifts},
 Ergodic Theory Dynam. Sytems
{\bf 29}(2009), pp.\ 1613--1639.
%

\bibitem{PRW2}
{\sc D. ~Pask, I.~Raeburn and N. A. ~Weaver},
{\it Periodic $2$-graphs arising from subshifts},
 preprint arXiv:0911.0730.







%






\bibitem{Pim}
{\sc M. V. Pimsner},
{\it A class of $C^*$-algebras generalizing both Cuntz--Krieger algebras 
and crossed products by ${\Bbb Z}$},
 in Free Probability Theory,  Fields Institute Communications
{\bf 12}(1996), pp. 189--212.


\bibitem{PCHR}
{\sc G. A. ~Pino, J. ~Clark, A. A. ~Huff and I.~Raeburn},
{\it Kumjian-Pask algebras of higher rank graphs},
 preprint arXiv:1106.4361v1.






 \bibitem{RoSte}
{\sc G. Robertson and Steger},
{\it  Affine buildings, tiling systems and higher rank Cuntz--Krieger algebras},
J.\ Reine Angew.\ Math.\
{\bf 513}(1999), pp.\ 115-144.



\bibitem{Ro}
{\sc M. R{\o}rdom},
{\it Classification of Cuntz-Krieger algebras},
 K-theory {\bf 9}(1995), pp.\  31--58.



















\bibitem{KSch}{\sc K. ~Schmidt},
{\it Dynamical systems of algebraic origin},
Birkh{\"a}user, Basel, Boston, Berlin 
(1995).

\bibitem{Wang}
{\sc H. ~Wang},
{\it Notes on a class of tiling problems},
Fundam. \ Math.\
{\bf 82}(1975), pp. 295--305.
 








\end{thebibliography}
\end{document}